\documentclass[a4paper,10pt]{paper}

\usepackage[utf8]{inputenc}
\usepackage{amsfonts,amsmath,amsthm,amssymb}
\usepackage[all]{xy}
\usepackage{verbatim}
\usepackage{mathrsfs}
\usepackage{yfonts}
\usepackage{a4wide}
\usepackage{color}
\usepackage{wallpaper}
\usepackage{fancyhdr}
\usepackage{subfig}

\title{Random walks on oriented lattices and Martin boundary}
\author{B. de Loynes}

\newtheorem{e-proposition}[theorem]{Proposition}

\newtheorem{e-definition}[theorem]{Definition\rm}

\newtheorem{theoreme}{Th\'eor\`eme}[section]

\newtheorem{proposition}[theoreme]{Proposition}

\newtheorem{definition}[theoreme]{D\'efinition\rm}

\begin{document}

\maketitle
\section*{Abridged English version}

In this note, we aim at sudying the Martin boundary of a simple random walk on a directed lattice. More precisely, let us denote by $(e_1,e_2)$ the canonical basis of $\mathbb{Z}^2$, and by $\mathsf{sgn}$ the sign function on the integers which is equal to $-1$ on negative integers and $1$ elsewhere. Thus, we are considering a $\mathbb{Z}^2$-valued Markov chain $(M_n)_{n \geq 0}$ with Markov operator $P$ defined for $(u,v) \in \mathbb{Z}^2 \times \mathbb{Z}^2$ as follows
\begin{enumerate}
\item $P(u,v) = \frac{1}{3}$ if $u_2 \neq 0$ and either $v=u \pm e_2$ or $v=u+\mathsf{sgn}(u_2) e_1$ ;
\item $P(u,v) = \frac{1}{2}$ if $v=u \pm e_2$ and $u_2 = 0$.
\end{enumerate}

For this random walk, the following can be shown.

\noindent
\textbf{Theorem 2.1}
\textit{The Martin boundary of the simple random walk $(M_n)_{n \geq 0}$ is trivial : all positive harmonic functions are constant.
}

In our case, the triviality of the Martin boundary comes from precise estimates of the Green function. Two kinds of behavior can be observed ; let $z \in \mathbb{Z} \times \{0\}$, then there exist positive constants $c$ and $s(\lambda)$ such that
\begin{enumerate}
\item $G(z,y) \sim s(\lambda) |y_2|^{-1}$ if $y_1y_2^{-2}$ converges to $\lambda \in \mathbb{R}$ ;
\item $G(z,y) \sim c |y_1-z|^{-\frac{1}{2}}$ if $y_1y_2^{-2}$ goes to $\pm \infty$.
\end{enumerate}
The fact that we give estimates of the Green function only for $z \in \mathbb{Z} \times \{0\}$ comes from technical difficulties which can be removed by the formula (2) of section \ref{originale}.

The triviality of the Martin boundary implies the triviality of the Poisson boundary, thus all bounded harmonic functions are constant (see \cite{kai} and references therein for a modern approach). Actually, ideas in \cite{cho} or in \cite{spi} can be suited to the Markov operator we are considering in this note and allows a direct proof of the triviality of the Poisson boundary. We do not follow this approach here because we obtain the stronger result of triviality of the Martin boundary.

\section{Introduction}
\label{}

Un \emph{di-graphe} $\mathbf{G}$ est la donn\'ee d'un couple $(\mathbf{V},\mathbf{E})$ d'un ensemble d\'enombrable $\mathbf{V}$ de points et d'un sous-ensemble $\mathbf{E} \subset \mathbf{V} \times \mathbf{V}$ d'ar\^etes orient\'ees. Si $e \in \mathbf{E}$ est une ar\^ete, on note $s(e)$ (resp. $t(e)$) le point \emph{source} (resp. le point \emph{terminal}). Plus pr\'ecis\'ement, $s(e)=u$ et $t(e)=v$ si $e=(u,v)$. Pour chaque point $u \in \mathbf{V}$, on d\'efinit le \emph{degr\'e sortant} de $u$, et on le note $d_u^-$ le nombre d'\'el\'ements de l'ensemble $\{ e \in \mathbf{E} : s(e)=u \}$. Un di-graphe est dit transitif si pour tout point $u,v \in \mathbf{V}$ il existe une suite finie de points $\{w_0, \cdots, w_k \}$ v\'erifiant $w_0=u$, $w_k=v$ et $(w_i,w_{i+1}) \in E$ pour tout $i \in \{0, \cdots, k-1\}$.
\begin{definition}
Soit $\mathbf{G}=(\mathbf{V},\mathbf{E})$ un di-graphe. La marche al\'eatoire simple sur $\mathbf{G}$ est la cha\^ine de Markov \`a valeurs dans $\mathbf{V}$, not\'ee $((M_n)_{n \geq 0},P,\mu)$ o\`u $\mu$ est une probabilit\'e sur $\mathbf{V}$, et  $P$ le noyau de Markov d\'efini par $P(u,v)=\frac{1}{d_u^-}$ si $(u,v) \in \mathbf{E}$ et $0$ sinon.
\end{definition}
Dans cette note, on consid\`ere les treillis bi-dimensionnels, \emph{i.e.} on supposera que $\mathbf{V}=\mathbb{Z}^2$ et que $\mathbf{E}$ est un sous-ensemble de l'ensemble des proches voisins dans $\mathbb{Z}^2$. On d\'ecompose $\mathbf{V}=\mathbf{V}_1 \times \mathbf{V}_2$ en les directions horizontales et verticales, c'est \`a dire, si $v \in \mathbf{V}$, on \'ecrit $v=(v_1,v_2)$ avec $v_1 \in \mathbf{V}_1$ et $v_2 \in \mathbf{V}_2$.
\begin{definition}
Soient $\mathbf{V}=\mathbf{V}_1 \times \mathbf{V}_2$ et $(\epsilon_y)_{y \in \mathbf{V}_2}$ \`a valeurs dans $\{-1,0,1\}$. On appelle treillis $\epsilon$-orient\'e, et on le note $\mathbf{G}=(\mathbf{G},\epsilon)$, le graphe orient\'e dont l'ensemble $\mathbf{V}$ de points est donn\'e par $\mathbf{Z}^2$ et $(u,v) \in \mathbf{E}$ est une ar\^ete orient\'ee si et seulement si l'une des conditions suivantes est v\'erifi\'ees :
\begin{itemize}
\item soit $v_1=u_1$ et $v_2=u_2 \pm 1$ ;
\item soit $v_2=u_2$ et $v_1 = u_1 + \epsilon_{u_2}$ si $\epsilon_{u_2} \neq 0$.
\end{itemize}
\end{definition}
Dans cette note, on \'etudiera le graphe $\mathbb{H}=(\mathbf{G},\epsilon)$ o\`u l'orientation $\epsilon$ est d\'efinie par $\epsilon_0=0$ et $\epsilon_y=\mathsf{sgn}(y)$ o\`u $\mathsf{sgn}$ est la fonction signe. Il est clair que $\mathbb{H}$ est transitif, ainsi la marche al\'eatoire simple sur $\mathbb{H}$ est irr\'eductible.

Les marches al\'eatoires sur ce type de graphe sont \'etroitement li\'ees aux diffusions en milieu poreux. On pourra lire avec int\'er\^et \cite{mat} o\`u ce point de vue est d\'evelopp\'e.


\section{R\'esultat principal}

Soit $(M_n)_{n \geq 0}$ la marche al\'eatoire simple sur $\mathbb{H}$. On notera $\mathbf{P}^\mu$ la probabilit\'e, sur l'espace des trajectoires, associ\'ee au noyau $P$ et plus simplement $\mathbf{P}^x$ si $\mu$ est la masse de Dirac en $x \in \mathbf{V}^2$. On d\'efinit par r\'ecurrence la suite de temps d'arr\^et $(\tau_n)_{n \geq 0}$ par $\tau_0=0$ et $\tau_{n+1}=\inf \{ t \geq \tau_{n}+1 : M_t^{(2)}=0 \}$ o\`u $M_t=(M_t^{(1)},M_t^{(2)})$. Il est facile de voir que $\tau_n < \infty$ $\mathbf{P}^x$-p.s. pour tout $n \geq 0$. La suite de variables al\'eatoires $(M_{\tau_n})_{n \geq 0}$ est alors elle-m\^eme une cha\^ine de Markov que l'on appellera \emph{cha\^ine de Markov induite}. Par d\'efinition, $M_{\tau_n}^{(2)}=0$, ainsi, si $\mu$ est \`a support dans $\mathbf{V}_1 \times \{ 0 \}$, not\'e par la suite $\mathcal{V}_0$, la cha\^ine de Markov induite est confin\'ee dans ce sous-ensemble $\mathcal{V}_0$. La cha\^ine de Markov induite sera alors, selon le contexte, consid\'er\'ee comme une cha\^ine de Markov \`a valeurs dans $\mathbb{Z}^2$ ou $\mathbb{Z}$.

Dans \cite{pet}, il est montr\'e que la cha\^ine de Markov originale est transiente, et cela est une cons\'equence de la transience de la cha\^ine de Markov induite qui peut-\^etre vue comme une marche al\'eatoire sur $\mathbb{Z}$ \`a sauts non born\'es. Ainsi, l'\'etude de la fronti\`ere de Martin de ces cha\^ines de Markov est pertinente.

\begin{theoreme}
Les fronti\`eres de Martin des cha\^ines de Markov induites et originales sont triviales, i.e. les seules fonctions harmoniques positives sont les fonctions constantes.
\end{theoreme}

\section{Fronti\`ere de Martin de la cha\^ine de Markov induite} \label{induite}

La compactification de Martin est longuement pr\'esent\'ee dans \cite{woe}. Essentiellement, la compactification de Martin fait intervenir le noyau de Martin : si $o$ est un point-base du graphe $\mathbb{H}$, le noyau de Martin est d\'efini pour $x,y \in \mathbb{H}$ par $K(x,y)=\frac{G(x,y)}{G(o,y)}$ o\`u $G$ est la fonction de Green associ\'ee \`a la cha\^ine de Markov. Ainsi, des estim\'ees pr\'ecises de $G$ permettent de donner une description de la fronti\`ere de Martin.

\begin{definition}
Soient $(\psi)_{n \geq 0}$ une suite de variables al\'eatoires ind\'ependantes et identiquement distribu\'ees \`a valeurs dans $\{-1,1\}$ de loi de Bernouilli de param\`etre $\frac{1}{2}$ et $(Y_n)_{n \geq 0}$ une suite de variables al\'eatoires d\'efinies par $Y_0=M_0^{(2)}$ et $Y_n=Y_0+\sum_{k=0}^{n-1} \psi_k$. On note $\eta_n(y)$ le temps local de $(Y_n)_{n \geq 0}$ en $y \in \mathbb{Z}$, i.e. $\eta_n(y)=\sum_{k=0}^{n} 1_{Y_k=y}$.
\end{definition}
\begin{definition}
Soit $(\sigma_n)_{n \geq 0}$ une suite de temps d'arr\^et d\'efinie par $\sigma_0=0$ et $\sigma_{n+1} = \inf \{ t \geq \sigma_n+1 : Y_t=0 \}$.
\end{definition}
\begin{definition}
Soit $(\xi_n^{(y)})_{n \geq 1, y \in \mathbf{V}_2}$ une suite de variables al\'eatoires ind\'ependantes et identiquement distribu\'ees \`a valeurs dans $\mathbb{N}$ de loi g\'eom\'etrique de param\`etre $p=\frac{1}{3}=1-q$ et pour $n \geq 0$ posons $X_n=\sum_{y \in \mathbf{V}_2} \epsilon_y \sum_{i=1}^{\eta_{\sigma_1-1}(y)} \xi_i^{(y)}$.
Enfin, on note par $T_n$ le temps d\'efini par $T_n=n+\sum_{y \in \mathbf{V}_2} \sum_{i=1}^{\eta_{\sigma_1-1}(y)} \xi_i^{(y)}$, avec la convention usuelle $\sum_\emptyset = 0$.
\end{definition}
Dans \cite{pet}, il est montr\'e que $M_{T_n}\stackrel{\textrm{\tiny loi}}{=}(X_n,Y_n)$ et $\tau_n\stackrel{\textrm{\tiny loi}}{=}T_{\sigma_n}$. Ceci permet de calculer la loi de $M_{\tau_n}$.
\begin{proposition}
La loi de $M_{\tau_n}$ est d\'etermin\'ee par sa fonction caract\'eristique $\mathbf{E}^0(e^{i \langle t,M_{\tau_n} \rangle })=\mathbf{E}^0(e^{it_1X_{\sigma_1}})^n$. De plus, la fonction caract\'eristique de $X_{\sigma_1}$ est donn\'ee par $\phi(t)= \mathsf{Re} \textrm{ } r(t)^{-1} g(r(t))$ o\`u $r(t)=[3-2e^{it}]^{-1}$ et $g(x)=\frac{1-\sqrt{1-x^2}}{x}$.
\end{proposition}
En particulier, la fonction de Green associ\'ee \`a la cha\^ine de Markov $(M_{\tau_n})_{n \geq 0}$ est donn\'ee, pour $u,v \in \mathcal{V}_0$, par $G_0(u,v)=G_0(0,v-u)=(2\pi)^{-1} \int_{-\pi}^\pi e^{-it(v-u)}[1-\phi(t)]^{-1} dt$. La fonction $[1-\phi(\cdot)]^{-1}$ poss\`ede, en $0$, une singularit\'e en $\frac{1}{\sqrt{|t|}}$. Une d\'ecomposition en s\'eries permet de traiter cette singularit\'e \`a part et de trouver l'\'equivalent suivant.
\begin{proposition} \label{k0}
Il existe une constante $c > 0$ telle que $G_0(u,v)=G_0(0,v-u) \sim c |v-u|^{-\frac{1}{2}}$ lorsque $|v|$ tend vers l'infini.
\end{proposition}
Ainsi le noyau de Martin $K_0(u,v)=\frac{G_0(u,v)}{G_0(0,v)}$ poss\`ede une unique valeur d'adh\'erence lorsque $|v|$ tend vers l'infini ce qui implique la trivialit\'e de la fronti\`ere de Martin.

\section{Fronti\`ere de Martin de la cha\^ine de Markov originale}
\label{originale}

Dans ce paragraphe on s'int\'eresse \`a la fronti\`ere de Martin de la cha\^ine originale. On notera $G$ (resp. $K$) la fonction de Green (resp. le noyau de Martin) associ\'ee \`a la marche originale.

On note $\eta_{s,t}(y)$, pour $s,t \geq 0$ et $y \in \mathbb{H}$, le temps local de $(M_n)_{n \geq 0}$ en $y$, \emph{i.e.} $\eta_{s,t}(y)=\sum_{k=s}^{t-1} 1_{M_n=y}$ toujours avec la convention $\sum_\emptyset=0$. On peut alors exprimer le noyau de Martin de la cha\^ine originale en fonction de celui de la cha\^ine induite.
\begin{displaymath}
K(x,y)= \left \{ \begin{array}{lr}
K_0(x,y) \textrm{ si } (x,y) \in \mathcal{V}_0 \textrm{, } & (1) \\
\frac{\mathbf{E}^x(\eta_{0,\tau_1}(y))}{G(0,y)} + \sum_{z \in \mathbb{Z}} \nu_x(z) K(z,y) \textrm{ sinon. } & (2) \\
\end{array} \right . 
\end{displaymath}
o\`u $\nu_x(z)=\mathbf{P}^x(M_{\tau_1}=(z,0))=\mathbf{P}^x(X_{\sigma_1}=z)$.

En utilisant les sym\'etries de $\mathbb{H}$, en particulier que ce graphe est invariant par les translations horizontales et par les sym\'etries centrales de centre situ\'e sur $\mathcal{V}_0$, on peut montrer la proposition suivante.
\begin{proposition}
Soient $z \in \mathcal{V}_0$ et $y \in \mathbb{H}$, alors la fonction de Green $G$ est donn\'ee par
\begin{displaymath}
G(z,y)=(2\pi)^{-1} \int_{-\pi}^\pi e^{it(y_1-z_1)} \frac{g(r(t))^{|y_2|}}{1-\phi(t)} dt,
\end{displaymath}
o\`u $r,g$ et $\phi$ sont d\'efinies dans le paragraphe \ref{induite}.
\end{proposition}
Il est alors possible d'obtenir les estim\'ees suivantes de la fonction de Green.
\begin{proposition} \label{greenestimates}
Soient $z \in \mathcal{V}_0$, $y \in \mathbb{H}$ et supposons que $y_1y_2^{-2}$ converge vers $\lambda \in \mathbb{R} \cup \{ \pm \infty \}$, alors il existe pour tout $\lambda$ une constante $s(\lambda)>0$ telle que
\begin{enumerate}
\item $G(z,y) \sim s(\pm \infty) |y_1-z_1|^{-\frac{1}{2}}$ lorsque $y_1y_2^{-2}$ tend vers $\pm \infty$ ;
\item $G(z,y) \sim s(\lambda) |y_2|^{-1}$ lorsque $y_1y_2^{-2}$ converge vers $\lambda$.
\end{enumerate}
\end{proposition}
Par un calcul direct, on peut montrer qu'il est possible de passer \`a la limite dans la somme $\sum_{z \in \mathbb{Z}} \nu_x(z) K(z,y)$, d'o\`u la proposition suivante.
\begin{proposition} \label{sumestimate}
Soit $x \in \mathbb{H}$, alors $\sum_{z \in \mathbb{Z}} \nu_x(z)K(z,y)$ poss\`ede une unique valeur d'adh\'erence lorsque $|y|$ tend vers l'infini.
\end{proposition}
\begin{proposition} \label{first}
La quantit\'e $\mathbf{E}^0(\eta_{0,\tau_1}(y))$ d\'ecro\^it vers $0$ en $o(|y_1|^{-\frac{1}{2}})$ lorsque $y_1y_2^{-2}$ tend vers $\pm \infty$ et en $o(|y_2|^{-1})$ lorsque $y_1y_2^{-2}$ converge.
\end{proposition}
Par les propositions \ref{greenestimates} et \ref{first}, le premier terme dans l'\'equation (2) converge vers $0$ lorsque $|y|$ tend vers l'infini ; et par la proposition \ref{sumestimate}, le deuxi\`eme terme de l'\'equation (2) converge vers $1$ lorsque $|y|$ tend vers l'infini. Enfin, le terme de l'\'equation (1) tend vers $1$ lorsque $|y|$ tend vers l'infini par la proposition \ref{k0}. Ainsi, le noyau de Martin $K$ poss\`ede une unique valeur d'adh\'erence \`a l'infini ce qui montre la trivialit\'e de la fronti\`ere de Martin de la cha\^ine originale.

Dans le cas o\`u $\epsilon$ est une suite de variables al\'eatoires ind\'ependantes et identiquement distribu\'ees, il est montr\'e dans \cite{pet}, que la marche al\'eatoire simple est encore transiente (pour presque toute orientation). Ce r\'esultat a \'et\'e g\'en\'eralis\'e dans \cite{gui} o\`u est consid\'er\'e le cas o\`u $(f_y)$ est une suite stationnaire \`a valeurs dans $[0,1]$ et o\`u $\epsilon_y$ est pris \'egal \`a 1 avec probabilit\'e $f_y$ et -1 avec probabilit\'e $1-f_y$. Sous la condition $\mathbf{E}(f_0(1-f_0))^{-1/2} < \infty$, la marche simple est transiente. Enfin, dans \cite{pen}, est consid\'er\'e le cas o\`u l'orientation $\epsilon$ est une suite stationnaire v\'erifiant certaines conditions de d\'ecorr\'elations. En moyenne, une orientation al\'eatoire est moins d\'eg\'en\'er\'ee que l'orientation choisie dans cette note, on peut donc s'attendre, dans ce cadre, \`a la trivialit\'e des fronti\`eres de Martin pour presque toute orientation.


\end{document}